%% file: rbfsi2.tex
\documentclass{amsart}

\usepackage{graphicx}
\usepackage{amsmath}
\usepackage{mathrsfs}
\usepackage{amssymb}
\usepackage{amsfonts}
\usepackage{algorithm}
\usepackage{algorithmic}
\usepackage{subfigure}

\newcommand{\param}{\boldsymbol{\mu}}

\newcommand{\bx}{\boldsymbol{x}}

\newcommand{\mup}{\boldsymbol{\mu}_{\ell,m}}
\newcommand{\Pp}{\boldsymbol{P}_{\ell,m}}
\newcommand{\Fp}{T}

\newcommand{\velospace}{\mathcal{Y}}
\newcommand{\prespace}{\mathcal{Q}}

\newcommand{\bm}[1] {\mbox{\boldmath $#1$\unboldmath}}

\begin{document}
\title[Reduced formulation of steady fluid-structure interaction]{Reduced formulation of a steady fluid-structure \\interaction problem with parametric coupling}

\author[Toni Lassila]{Toni Lassila${}^*$}

\author[Gianluigi Rozza]{Gianluigi Rozza${}^{\dagger}$}

\thanks{${}^*$ Department of Mathematics and Systems Analysis, Aalto University School of Science and Technology. email: \texttt{toni.lassila@tkk.fi}
\\ \\
${}^{\dagger}$ Modelling and Scientific Computing Chair, Mathematics Institute of Computational Science and Engineering, {\'{E}}cole Polytechnique F{\'{e}}d{\'{e}}rale de Lausanne. email: \texttt{gianluigi.rozza@epfl.ch}
}

\begin{abstract} 
 We propose a two-fold approach to model reduction of fluid-structure interaction.
  The state equations for the fluid are solved with reduced basis methods. 
  These are model reduction methods for parametric partial
  differential equations using well-chosen snapshot solutions in order to build a set of global
  basis functions. The other reduction is in terms of the geometric complexity of the
  moving fluid-structure interface. We use free-form deformations to parameterize the
  perturbation of the flow channel at rest configuration. As a computational example we
  consider a steady fluid-structure interaction problem: an incmpressible Stokes flow in a
  channel that has a flexible wall.
\end{abstract}

\maketitle

\section{INTRODUCTION}
 
  Many problems in fields such as aerodynamics and biomechanics can be expressed
  as fluid-structure interaction (FSI) problems. The mathematical modelling of such coupled
  problems consists of four main parts: solution of the fluid equations given the current 
  fluid geometry, solution of the structural displacements given the normal stresses exerted 
  by the fluid, fulfillment of coupling constraints to achieve force balance across the 
  interface, and transport of the fluid-structure interface. Even when both the fluid and 
  structure equations are independently linear, the geometric variability of the fluid-structure 
  interface yields a significant nonlinearity in the coupled system. Solution methods for 
  fluid-structure interaction problems are iterative in nature and involve the repeated 
  solution of the fluid and structure equations in many different configurations.

  Fluid-structure interaction problems arise in modelling of the arterial deformations in
  the human cardiovascular system. A comprehensive treatment of FSI strategies in
  cardiovascular modelling can be found in \cite{NobileThesis,QuarteroniFormaggia}. Blood
  flow is pulsatile and the displacements of the arterial walls are relatively large. This
  leads to stability considerations that necessitate the use of implicit solver strategies
  of the coupled FS problem, again increasing computational complexity. 
  
  Model Order Reduction (MOR) is a cross-disciplinary field that strives to systematically
  reduce the dynamics of a complicated ordinary or partial differential equation model to
  a simpler and computationally more tractable one. Model reduction techniques have also 
  been proposed for nonlinear systems. One particularly common approach is the 
  Proper Orthogonal Decomposition (POD) technique that aims at decomposing the dynamics of
  a time-dependent system into fundamental modes, and then choosing only very few of the
  most important modes to represent the entire dynamics of the system. 

  A modal reduction of the structural equations was adopted in \cite{Murea} to reduce
  the complexity of an FSI problem involving pulsatile flow in a channel, but no attempt to reduce
  the complexity of solving the fluid equations was made. In another work by the same author
  \cite{Murea2} the steady problem was parameterized using shape functions for the
  boundary trace of the pressure obtained from an assumption of near-Poiseuille flow profile.
  The strong coupling of normal stresses on the interface was treated in the reduced space
  obtained by considering only the leading structural eigenmodes, and solved as a 
  least-squares optimization problem.
  
  Another approach consists of the Reduced Basis (RB) methods, which is generally used on 
  parametrized partial (and ordinary) differential equations. RB methods are based on a 
  greedy sampling algorithm that chooses snapshot solutions of the parametric PDE at 
  different parameter values for constructing a global approximation basis. For an introduction, see
  \cite{PateraRozza} and \cite{RozzaEtAl}. The reduction of steady incompressible Stokes
  flow (without FSI) in an arterial bypass configuration by parameterization with piecewise
  affine maps was considered in \cite{RozzaVeroy}. The same
  problem parameterized with nonaffine maps was considered in \cite{Rozza} for Stokes
  equations and in \cite{QuarteroniRozza} for Navier-Stokes equations. We also
  mention the reduced basis element method proposed for Stokes flows in \cite{LovgrenEtAl}
  that was derived to address complex flow networks with multiple branches and
  junctions. A posteriori error estimation for Navier-Stokes equations solved by reduced
  basis methods is considered in \cite{Deparis,VeroyPatera}. We use reduced basis methods
  to approximate the solution of the Stokes equations in a parametric flow geometry.

\section{PROBLEM OF STEADY STOKES FLOW IN A FLEXIBLE CHANNEL} \label{sec:fsi}

\begin{figure}
    \centering
    \input{fsi.pstex_t}
    \caption{Stokes flow in a channel with a flexible wall section}
    \label{fig:channel}
  \end{figure}
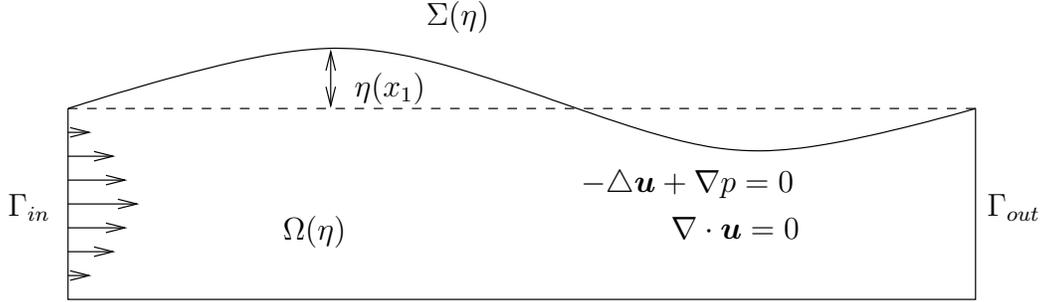

  We denote by $H^m(X)$ the usual Sobolev space of real-valued functions on $X$ with $m$th 
  weak derivatives being square-integrable. Let $\Omega_0 \subset \mathbb{R}^2$ be a bounded 
  domain with Lipschitz boundary $\partial \Omega_0$ that represents the fluid domain at 
  rest configuration (without deformations induced by the fluid flow). We denote the flexible 
  part of the boundary as $\Sigma_0 \subset \partial \Omega_0$. We assume the geometry 
  displayed in Fig. \ref{fig:channel} -- a straight 2-d tube, where the displacement of the
  upper wall is described with a function
  $\eta \in H^1_0(\Sigma_0) := \{ \phi \in H^1(\Sigma_0) \: : \: \phi(a) = \phi(b) = 0 \}$.
  The tube is at ``rest" configuration when $\eta \equiv 0$. We denote by $\Omega(\eta)$ the
  deformed domain. In the deformed domain we have the incompressible steady Stokes problem to 
  find $\bm u \in [H^1(\Omega(\eta))]^2$ and $p \in L^2_0(\Omega(\eta))$ s.t.
  \begin{equation} \label{eq:stokes}
  \left\{
    \begin{aligned}
    \int_{\Omega(\eta)} \left[
     \nu \nabla \bm u \cdot \nabla \bm v - p \nabla \cdot \bm v \right] \: d\Omega &= \int_{\Omega(\eta)} \bm f^F \cdot \bm v, \quad &\textrm{ for all } \bm v \in [H^1(\Omega(\eta))]^2 \\
    \int_{\Omega(\eta)} q \nabla \cdot \bm u \: d\Omega &= 0, \quad &\textrm{ for all } q \in L^2_0(\Omega(\eta)) \\
    \bm u = \bm u_0 \textrm{ on } \partial \Omega(\eta) \: \backslash \: \Sigma(\eta), &\quad \bm u = \bm 0 &\textrm{ on } \Sigma(\eta)
    \end{aligned}    
    \right. 
  \end{equation}
  where $\nu$ is the viscosity, and $\bm f^F$ is the volume force acting on the fluid. 
  The Dirichlet data are assumed to satisfy a conservation principle
  \begin{equation} \label{eq:conservation}
    \int_{\partial \Omega(\eta)} \bm u_0 \cdot \bm n \: d\Gamma = 0.
  \end{equation}

  The fluid equations are then to be coupled together with the structural equations given the coupling
  condition that the displacement of the structure is the same as would be obtained from the structural
  equations, where the RHS is the stress exerted by the fluid on the interface. 
  The structural equations can be chosen in several different ways (in \cite{Grandmont} and \cite{Murea}
  different fourth order models are considered). We choose, after a suitable assumption of small 
  structural displacement, to model the displacement of the wall using a 1-d elliptic equation 
  (``elastic membrane") to find $\eta \in H^1_0(\Sigma_0)$ s.t.
  \begin{equation}
    \int_{\Sigma_0} K(\bm x) \eta' \phi' \: d\Gamma = \int_{\Sigma_0} \tau(\bm u, p) \phi \: d\Gamma, \quad \forall \: \phi \in H^1_0(\Sigma_0)
  \end{equation}
  where $\tau(\bm u, p)$ is the normal traction on the interface imposed by the fluid
  \begin{equation}
    \tau(\bm u, p) = \left[ p \bm n - \nu \left( \nabla \bm u + \nabla \bm u^T \right) \cdot \bm n \right]^T \begin{bmatrix} 0 \\ 1 \end{bmatrix}
  \end{equation}
  obtained from the solution $(\bm u, p)$ of the fluid equation (\ref{eq:stokes}), and 
  $K(\bm x) \geq K_0 > 0$ is the spring constant of the wall to displacements in the normal direction. 
 
  In operator form, the coupled fluid-structure problem can be considered to be of the form
  \begin{equation} \label{eq:system}
    \begin{bmatrix} -\nu \triangle & \nabla & 0 \\ \textrm{ div} & 0 & 0 \\ 0 & 0 & K \partial_{xx} \end{bmatrix}
    \begin{bmatrix} \bm u \\ p \\ \eta \end{bmatrix} +
    \begin{bmatrix} 0 \\ 0 \\ \tau(\bm u, p) \end{bmatrix}
    = 
    \begin{bmatrix} \bm f^F \\ 0 \\ 0 \end{bmatrix}.
  \end{equation}   
  Because the fluid domain depends on the domain defined by the function $\eta$, this is
  in fact a nonlinear free-boundary problem. 

\section{PARAMETERIZATION WITH FREE-FORM DEFORMATIONS} \label{sect:ffd}

  To reduce the complexity of the free-boundary problem (\ref{eq:system})
  we introduce a parameterization of the fluid domain, $\Omega(\param)$, where the
  parameters $\param \in \mathcal{D}$ belong to some low-dimensional parameter space.
  One option to construct $\Omega(\param)$ would be to use a family of parametric 
  curves, such as B-splines, to directly parameterize the displacement function $\eta$.
  We prefer to use \emph{free-form deformations} \cite{SederbergParry}.
  
\begin{figure}
    \centering
    \includegraphics[width=10cm]{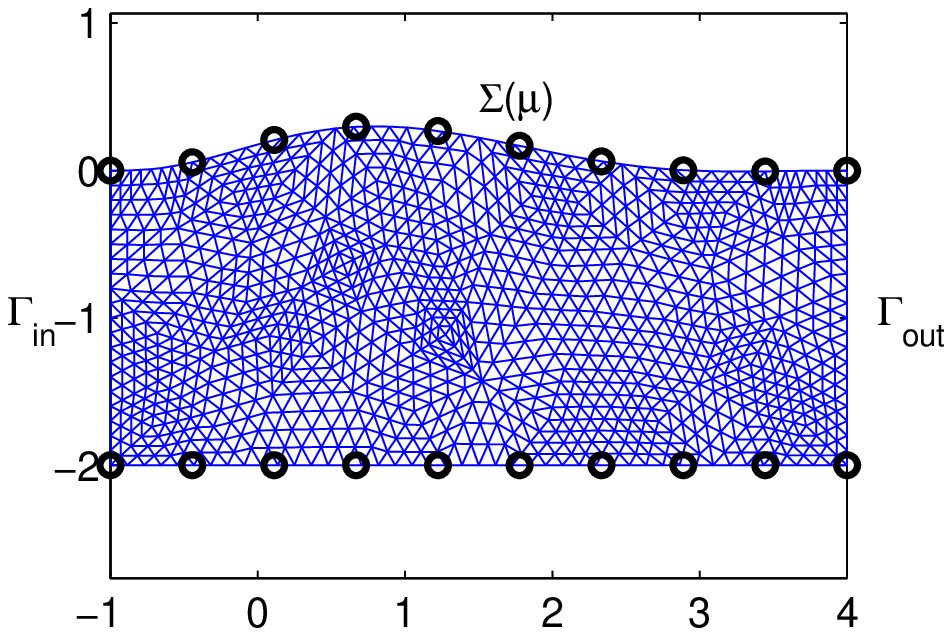}
    \caption{Free-form deformation of the flexible channel with six parameters.}
    \label{fig:ffdchannel}
  \end{figure}

  Free-form deformations are a technique of obtaining small, parametric deformations of an
  arbitrary reference domain $\Omega_0 \subset \mathbb{R}^d$ for $d=2,3$. First we 
  map the reference domain to a subset of the unit square with a continuous invertible map 
  $\Psi : \Omega_0 \to (0,1) \times (0,1)$ and then overlay on $\Psi(\Omega_0)$ 
  a regular grid of control points $\Pp^0$, $\ell=0,\ldots,L$ and $m=0,\ldots,M$, so that
  $\Pp^0 = \begin{bmatrix} \ell / L & m / M \end{bmatrix}^T$.
  The parameterization is obtained by allowing a subset of the control points to move and then using a
  spline basis to construct a smooth deformation map as a function of the positions of the control
  points. The perturbed control points are given by a set of $(L+1)(M+1)$
  parameter vectors $\mup$. For each $\param$, a parametric domain map
  $\Fp : \Omega_0 \times \mathcal{D} \to \Omega(\param)$ is defined as
  \begin{equation} \label{eq:ffd}
    \Fp(\bx;\param) = \Psi^{-1} \left( \sum_{\ell=0}^L \sum_{m=0}^M b_{\ell,m}^{L,M}(\Psi(\bx)) \left( \Pp + \mup \right) \right),
  \end{equation}
  where $b_{\ell,m}^{L,M}$ are suitable polynomial basis functions. Here we take
  \begin{equation*}
    b_{\ell,m}^{L,M}(s,t) = \binom{L}{\ell}\binom{M}{m} (1-s)^{(L-\ell)} s^{\ell} (1-t)^{(M-m)} t^m,
  \end{equation*}
  that is the tensor products of the 1-d Bernstein basis polynomials 
  \begin{equation*}
    b_{\ell}^L(s) = \binom{L}{\ell} (1-s)^{(L-\ell)} s^{\ell}, \quad b_{m}^M(t) = \binom{M}{m} (1-t)^{(M-m)} t^m
  \end{equation*}
  defined on the unit square with local variables $(s,t) \in (0,1) \times (0,1)$.  
  It is typical that only some of the control points are allowed to move freely, we do not use the
  full parameterization but only a few of the most relevant shape parameters. 
  To parameterize the flow channel we use a $10 \times 2$ grid of control points.
  Six control points in the top row are moving in the $x_2$-direction, resulting in a
  smooth parameterized deformation of the channel with 6 real parameters.
  The result when applied on a sample computational mesh is shown in Fig. \ref{fig:ffdchannel}.
  
  After the parametric domain $\Omega(\param)$ is obtained as the image of a fixed reference domain
  $\Omega_0$ under the map $T : \Omega_0 \times \mathcal{D} \to \Omega(\param)$, which is assumed 
  to be differentiable and invertible, we proceed to transform the PDE on a parametric domain to a 
  parametric PDE on the fixed reference domain. Denoting by $J_{T_{\eta}} := \nabla_{\bx} T_{\eta}$ 
  the Jacobian matrix of $T_{\eta}$ we define the transformation tensors for the viscous term
  \begin{equation}
	\nu_{T_{\eta}} := J_{T_{\eta}}^{-T} J_{T_{\eta}}^{-1} \det(J_{T_{\eta}})
  \end{equation}
  and the pressure-divergence term
  \begin{equation}
	\chi_{T_{\eta}} := J_{T_{\eta}}^{-1} \det(J_{T_{\eta}})
  \end{equation}
  respectively.
  The inhomogeneous Dirichlet conditions are handled by lifting to a space of proper boundary conditions
  the solution of the homogeneous Stokes problem on the reference domain to find 
  $\bm u(\eta) \in [H^1_0(\Omega_0)]^2$ and $p(\eta) \in L^2_0(\Omega_0)$ s.t.
  \begin{equation} \label{eq:param_stokes}
  \left\{
    \begin{aligned}      
      \int_{\Omega_0} 
      \left( \nu
      \frac{\partial u_k}{\partial x_i} [\nu_{T_{\eta}}]_{i,j} \frac{\partial v_k}{\partial x_j}
       + p \: [\chi_{T_{\eta}}]_{k,j} \frac{\partial v_k}{\partial x_j} 
       \right) 
       \: d\Omega_0
       &= \int_{\Omega_0} \det(J_{T_{\eta}}) [f^F + f_{\textrm{lift}}]_k \: d\Omega_0, \\
       & \quad \quad \forall \: \bm v \in H^1_0(\Omega_0) \times H^1_0(\Omega_0) \\
      \int_{\Omega_0} 
      q \: [\chi_{T_{\eta}}]_{k,j} \frac{\partial u_k}{\partial x_j}       
       \: d\Omega_0 &= 0, \\
       & \quad \quad \forall \: q \in L^2_0(\Omega_0)
    \end{aligned}    
    \right.
  \end{equation}
  with summation understood over the indices $i=1,2$ and $j=1,2$. 
     
\section{PARAMETRIC FLUID-STRUCTURE COUPLING ALGORITHM} \label{sec:coupling}
  
  We can alternatively formulate the fluid-structure problem (\ref{eq:system}) as a
  calculus of variations problem to find $\eta \in H^{1}_0(\Sigma_0)$ the minimizer of
  \begin{equation} \label{eq:minproblem}
    \min_{\eta} \quad \tfrac{1}{2} \int_{\Sigma_0} | K \partial_{xx} \eta + \tau(\bm u, p) |^2 \: d\Gamma,\\
  \end{equation}
  \begin{equation*}
    \textrm{s.t.} \quad 
    \left\{
    \begin{aligned}
      \quad -\nu \triangle \bm u + \nabla p &= \bm f^F \quad &\textrm{ in } \Omega(\eta) \\
        \nabla \cdot \bm u &= 0 \quad &\textrm{ in } \Omega(\eta) 
    \end{aligned}
    \right.
  \end{equation*}
  If a solution of (\ref{eq:system}) exists then it is also a minimizer of (\ref{eq:minproblem}). 
  Problem (\ref{eq:minproblem}) in the parametric form given by (\ref{eq:param_stokes})
  is still valid: find $\param \in \mathcal{D}$ s.t.
  \begin{equation} \label{eq:minproblem2}
    \min_{\param} \quad \tfrac{1}{2} \int_{\Sigma_0} | K \partial_{xx} \eta + \tau(\bm u, p) |^2 \: d\Gamma,\\
  \end{equation}
  \begin{equation*}
    \textrm{s.t.} \quad 
    \left\{
    \begin{aligned}
      \quad -\nu \triangle \bm u + \nabla p &= \bm f^F \quad &\textrm{ in } \Omega(\param) \\
      \nabla \cdot \bm u &= 0 \quad &\textrm{ in } \Omega(\param) 
    \end{aligned}
    \right.
  \end{equation*}  
  but this time we expect that the optimal value of the functional is $J^* > 0$ and the coupling
  of stresses on the interface between fluid and structure is only achieved in a least-squares sense.
  The ``goodness of fit" depends on the dimension of the parameter space $\mathcal{D}$.
  We call this the \emph{strong parametric coupling} approach.

  It should be noted that we perform a reduction by parameterizing the displacement
  $\eta$ of the fluid-structure interface. In \cite{Murea2} the traction $\tau$ (which
  corresponds to the boundary trace of the pressure in this case) was instead
  parameterized. This alternative approach allowed the author to solve directly a
  parametric minimization problem of the form (\ref{eq:minproblem2}). Instead, we must
  consider the minimization problem
  \begin{equation} \label{eq:minproblem3}
    \min_{\param} \quad \tfrac{1}{2} \int_{\Sigma_0} \left[ | \eta(\param) - \hat{\eta} |^2 + | \partial_x \eta(\param) - \partial_x \hat{\eta} |^2 \right]  \: d\Gamma,\\
  \end{equation}
  \begin{equation*}
    \textrm{s.t.} \quad 
    \left\{
    \begin{aligned}
      \quad -\nu \triangle \bm u + \nabla p &= \bm f^F \quad \textrm{ in } \Omega(\param) \\
      \nabla \cdot \bm u &= 0 \quad \textrm{ in } \Omega(\param) \\
      K \partial_{xx} \hat{\eta} + \tau(\bm u, p) &= 0 \quad \textrm{ on } \Sigma(\param) 
    \end{aligned}
    \right..
  \end{equation*}  
  The reason for
  this is that, in order to reduce the complexity of the fluid equations, we must be able
  to transfer the parametric dependence to the Stokes equations. This can be done
  more easily by parameterizing the displacement rather than the boundary trace of the
  pressure.
  
  A computational algorithm for the solution of (\ref{eq:minproblem3}) is detailed in
  Algorithm 1.
  \begin{algorithm}[!hpbt]
	\caption{Strong parametric coupling of fluid and structure} 
	\label{algo:strongcoupling}
	\begin{algorithmic}[1]
		\REQUIRE initial guess $\param^0$
		\STATE Let $k=0$.
		\REPEAT
		\STATE Solve the discretized fluid equations for $(\bm u^k, p^k)$ in $\Omega(\param^k)$.
		\STATE Compute assumed traction $ \hat{\tau} = p^k|_{\Sigma(\param)} - \nu\left[ \nabla \bm u^k + [\nabla \bm u^k]^T \right]$ from $(\bm u^k,p^k)$.	 
        \STATE Solve the minimization problem 
        \[
		\param^{k+1} := \mathop{\textrm{argmin}}_{\param \in \mathcal{D}} \quad \tfrac{1}{2} \int_{\Sigma_0} \left\{ |\eta(\param) - \hat{\eta}|^2 + |\partial_x \eta(\param) - \partial_x \hat{\eta}|^2 \right\} \: d\Gamma
        \]
        to obtain a new configuration parameter,
      where $\eta(\param)$ is the displacement of the interface given by the parameterization
      of the geometry and the \emph{assumed displacement} $\hat{\eta} \in H^1_0(\Sigma_0)$ is the solution of
		  \[
    \displaystyle \int_{\Sigma_0} K \hat{\eta}' \phi' \: d\Gamma = \displaystyle \int_{\Sigma_0} \hat{\tau} \phi \: d\Gamma \quad \forall \: \phi \in H^1_0(\Sigma_0)
      \]
		\STATE Set $k := k + 1$.
		\UNTIL{stopping criteria $|\param^{k+1} - \param^k| < \varepsilon$ is met}
	\end{algorithmic}
  \end{algorithm}
  This is a fixed point algorithm in the parameter space. The algorithm's computationally most 
  intensive part is the solution of the parameterized fluid equations. We next discuss the reduced 
  basis method for approximation of the fluid solution $(\bm u^k, p^k)$.

\section{REDUCED BASIS APPROXIMATION OF THE STOKES PROBLEM}

  The computational efficiency of solving the coupling problem 
  (\ref{eq:minproblem3}) hinges on the efficient solution of the fluid equations on the
  parametric domain $\Omega(\param)$. We present briefly the general approach to reduced
  basis approximations for the steady Stokes problem. For more details see
  \cite{Rozza} and \cite{RozzaVeroy}.

  We denote by $\velospace := H^1(\Omega_0) \times H^1(\Omega_0)$ the velocity space,
  $\prespace := L^2(\Omega_0)$ the pressure space, and the continuous parametric
  bilinear forms $\mathcal{A} : \velospace \times \velospace \times \mathcal{D} \to \mathbb{R}$
  and $\mathcal{B} : \prespace \times \velospace  \times \mathcal{D} \to \mathbb{R}$.
  The weak form (\ref{eq:param_stokes}) of the parametric Stokes equations for velocity 
  and pressure can be written: given any $\param \in \mathcal{D}$, find 
  $\bm u(\param) \in \velospace$ and $p(\param) \in \prespace$ s.t.
  \begin{equation} \label{eq:problem1}
    \begin{aligned}
      \mathcal{A}(\bm u,\bm v; \param) + \mathcal{B}(p,\bm v;\param) &= \langle \bm F(\param), \bm v \rangle_{\velospace} \quad &\textrm{ for all } \bm v \in \velospace \\
      \mathcal{B}(q,\bm u(\param)) &= \langle G(\param), q \rangle_{\prespace} \quad &\textrm{ for all } q \in \prespace,
  \end{aligned}    
\end{equation}
  where $\bm F(\param) \in \velospace' \times \mathcal{D}$ and $G(\param) \in \prespace'  \times \mathcal{D}$
  are continuous parametric linear forms. The inhomogeneous Dirichlet condition is treated by
  introducing an extension of the boundary data $\tilde{\bm u}_0 \in H^1(\Omega_0)$
  and looking for a solution $\tilde{\bm u} \in H^1_0(\Omega_0) \times H^1_0(\Omega_0) \times \mathcal{D}$  
  to (\ref{eq:problem1}) with a modified right-hand side
  $\langle \tilde{\bm F}(\param), \bm v \rangle_{\velospace} := \langle \bm F(\param), \bm v \rangle_{\velospace} - \langle \mathcal{A}(\param) \tilde{\bm u}_0, \bm v \rangle_{\velospace}$. The solution is then
  recovered as $\bm u(\param) = \tilde{\bm u}(\param) +  \tilde{\bm u}_0$.

  Define $\alpha(\param)$ the parameter-dependent coercivity constant of the problem as
  \begin{equation}
    \alpha(\param) := \inf_{\bm v \in \velospace} \frac{\mathcal{A}(\bm v,\bm v;\param)}{||\bm v||^2_{\velospace}},
  \end{equation}
  and $\beta(\param)$ the parameter-dependent inf-sup constant of the problem as
  \begin{equation}
    \beta(\param) := \inf_{q \in \prespace} \sup_{\bm v \in \velospace} \frac{\mathcal{B}(q,\bm v;\param)}{||q||_{\prespace} ||\bm v||_{\velospace}}.
  \end{equation}
  A necessary condition known as the Ladyzhenskaya-Babu{\v{s}}ka-Brezzi condition for
  problem (\ref{eq:problem1}) to have a unique solution is that $\alpha(\param) \geq \alpha_0 > 0$ and
  $\beta(\param) \geq \beta_0 > 0$. We assume in what follows that the coercivity constant is always positive
  and concentrate on the inf-sup constant.

  From here on we forego any interest in the infinite-dimensional problem and assume instead 
  that both $\velospace$ and $\prespace$ are replaced by strictly finite-dimensional spaces 
  $\velospace_h$ and $\prespace_h$. More concretely, we consider only the finite-dimensional 
  problem obtained by discretization of the partial differential equations using the finite 
  element method. The discretized model problem
  \begin{equation} \label{eq:problem3}
    \begin{aligned}
      \mathcal{A}(\bm u_h, \bm v_h; \param) + \mathcal{B}(p_h, \bm v_h; \param) &= \langle \bm F_h(\param), \bm v_h \rangle_{\velospace_h} \quad &\textrm{ for all } \bm v_h \in \velospace_h \\
      \mathcal{B}(q_h, \bm u_h; \param) &= \langle G_h(\param), q_h \rangle_{\prespace_h} \quad &\textrm{ for all } q_h \in \prespace_h,
    \end{aligned}    
  \end{equation}
  involves a similar inf-sup condition to guarantee the existence of a unique
  velocity-pressure pair satisfying (\ref{eq:problem3}):
  \begin{equation} \label{eq:lbbcondition}
    \beta_h(\param) = \inf_{q_h \in \prespace_h} \sup_{\bm v_h \in \velospace_h} \frac{\mathcal{B}(q_h,\bm v_h; \param)}{||q_h||_{\prespace} ||\bm v_h||_{\velospace}} \geq \beta_0 > 0.
  \end{equation}
  In the case of the steady Stokes problem we employ the $\mathbb{P}^2-\mathbb{P}^1$
  Taylor-Hood elements \cite{GiraultRaviart,QuarteroniValli} that satisfy the condition
  (\ref{eq:lbbcondition}) without the necessity of adding extra stabilizing terms.

  The reduced basis method seeks an approximation of the finite element solution to
  (\ref{eq:problem3}). The first step is to construct a set of global basis functions.
   Let $\param^1,\ldots,\param^N$ be a set of snapshot parameter values chosen
  according to some rule. We denote by $\bm u(\param^n)$ and $p(\param^n)$ the
  corresponding snapshot solutions for the velocity and pressure. The snapshot
  solutions are obtained as solutions of a stable finite element formulation of
  (\ref{eq:problem3}). The reduced basis velocity space is then defined as 
  \begin{equation}
    \textrm{span}(u(\param^1),\ldots,u(\param^N)) =: \velospace_N \subset \velospace_h,
  \end{equation}
  and the reduced basis pressure space as 
  \begin{equation}
    \textrm{span}(p(\param^1),\ldots,p(\param^N)) =: \prespace_N \subset \prespace_h.
  \end{equation}
  We define the \emph{reduced basis Galerkin} problem: for any $\param \in \mathcal{D}$
  find $u_N \in \velospace_N$ and $p_N \in \prespace_N$ such that
  \begin{equation} \label{eq:problem4}
    \begin{aligned}
      \mathcal{A}(u_N(\param), v; \param) + \mathcal{B}(p_N(\param), v; \param) &= \langle \bm F_h(\param), v \rangle \quad &\textrm{ for all } v \in \velospace_N \\
      \mathcal{B}(q, u_N(\param); \param) &= \langle G_h(\param), q \rangle \quad &\textrm{ for all } q \in \prespace_N.
    \end{aligned}    
  \end{equation}
  Then $(\bm u_N(\param), p_N(\param))$ is a reduced basis approximation for 
  $(\bm u_h(\param), p_h(\param))$ the solution of (\ref{eq:problem3}). To represent and
  solve these equations in matrix form, an orthonormal basis is constructed for
  $\velospace_N$ and $\prespace_N$ to guarantee algebraic stability as $N \to \infty$ \cite{RozzaVeroy}.

  Provided that the dimension $N$ is chosen to be much smaller than the dimension of the
  finite element spaces, this system is inexpensive to solve. If the dependence of the bilinear
  forms on the parameters is smooth and the snapshot parameter values $\param_n$ are properly chosen, 
  we find that choosing a moderately small $N$ permits rapid convergence towards the finite element solution.
  For details on how to choose the snapshots with the help of a greedy algorithm that controls the 
  residual a posteriori approximation error, we refer to existing works on reduced basis methods
  \cite{GreplPatera,PateraRozza,RozzaEtAl}.

  A key issue in the reduced basis approximation of the Stokes problem is approximation
  stability. We remark that use of inf-sup-stable finite elements for the snapshot
  computation does not necessarily result in a stable reduced basis approximation when the
  basis functions are orthonormalized, as was observed in \cite{RozzaVeroy}. On the other 
  hand, the orthonormalization is vital for the algebraic stability of the reduced basis
  approximation. Thus there is a need to extend the concept of inf-sup constants
  into the reduced basis formulation. This is typically done by defining the so-called
  supremizers in velocity space and enriching the reduced basis velocity space
  $\velospace_N$ accordingly. Following \cite{RozzaVeroy} we define the inner supremizer
  operator $T^{\param} : \prespace_h \to \velospace_h$ as the solution of the elliptic
  problem
  \begin{equation}
    \langle T^{\param} q, \bm v \rangle_{\velospace} = \mathcal{B}(q, \bm v; \param) \quad \textrm{ for all } \bm v_h \in \velospace_h.
  \end{equation}
  After the pressure approximation space $\prespace_N$ and velocity approximation space
  $\velospace_N$ have been constructed we add all the velocity supremizers, i.e. velocity
  fields of the form $T^{\param} (p(\param^n))$, to the velocity space:
  \begin{equation}
    \tilde{\velospace}_N := \velospace_N \: \cup \: \textrm{span}\left[ T^{\param} (p(\param^1)),\ldots,T^{\param} (p(\param^N)) \right],
  \end{equation}
  where $\tilde{\velospace}_N$ is a stabilized reduced velocity space. Then
  \begin{equation}
    \tilde{\beta}_N(\param) :=  \inf_{q \in \prespace_N} \sup_{\bm v \in \tilde{\velospace}_N} \frac{\mathcal{B}(q, \bm v; \param)}{||q||_{\prespace_N} ||\bm v||_{\tilde{\velospace}_N}} \geq \beta_h(\param) \geq \beta_0 > 0.
  \end{equation}
  As a consequence of this approach to stabilization the dimension of the velocity
  approximation space grows to $2 N$, where $N$ is the original dimension of velocity and
  pressure approximation spaces respectively.

  The computational benefits of the reduced basis method rely on the assumption that the
  parametric bilinear forms are affinely parameterized:
  \begin{equation} \label{eq:affine}
      \mathcal{A}(\bm u, \bm v; \param) = \sum_{m=1}^{M_a} \Theta_a^m(\param) \mathcal{A}^m(\bm u, \bm v), \quad
      \mathcal{B}(p, \bm v; \param) = \sum_{m=1}^{M_b} \Theta_b^m(\param) \mathcal{B}^m(p, \bm v) 
  \end{equation}
  for some parameter-independent continuous bilinear forms $\mathcal{A}^m$ and $\mathcal{B}^m$,
  and parametric coefficient functions $\Theta_a^m$ and $\Theta_b^m$ that can be readily evaluated.
  Then the finite element equations (\ref{eq:problem3}) split into
  \begin{equation} \label{eq:problem7}
    \begin{aligned}
      \sum_{m=1}^{M_a} \Theta_a^m(\param) \mathcal{A}^m(\bm u_N(\param), \bm v) 
      + 
      \sum_{m=1}^{M_b} \Theta_b^m(\param) \mathcal{B}^m(p_N(\param), \bm v) 
	  &= \langle \bm F_h(\param), \bm v \rangle \quad &\textrm{ for all } \bm v \in \velospace_h \\
	\sum_{m=1}^{M_b} \Theta_b^m(\param) \mathcal{B}^m(q, \bm u_N(\param)) 
    &= \langle G_h(\param), q \rangle \quad &\textrm{ for all } q \in \prespace_h,
    \end{aligned}    
  \end{equation}
  and the reduced basis equations (\ref{eq:problem4}) respectively into
  \begin{equation} \label{eq:problem5}
    \begin{aligned}
      \sum_{m=1}^{M_a} \Theta_a^m(\param) \mathcal{A}^m(\bm u_N(\param), \bm v) 
      + 
      \sum_{m=1}^{M_b} \Theta_b^m(\param) \mathcal{B}^m(p_N(\param), \bm v) 
	  &= \langle \bm F_h(\param), \bm v \rangle \quad &\textrm{ for all } \bm v \in \velospace_N \\
	\sum_{m=1}^{M_b} \Theta_b^m(\param) \mathcal{B}^m(q, \bm u_N(\param)) 
    &= \langle G_h(\param), q \rangle \quad &\textrm{ for all } q \in \prespace_N.
    \end{aligned}    
  \end{equation}
  Once bases $\{ \bm \xi_n^v \}_{n=1}^{2N}$ and $\{ \xi_n^p \}_{n=1}^{N}$ for the
  reduced basis spaces $\tilde{\velospace}_N$ and $\prespace_N$ respectively have been
  constructed, the reduced system matrices
  \begin{equation}
    \begin{aligned}\,
    [A^m]_{k,\ell} &:= \mathcal{A}^m(\bm \xi_k^v,\bm \xi_{\ell}^v), \quad 1 \leq k, \ell \leq 2N; \\
    [B^m]_{k,\ell} &:= \mathcal{B}^m(\bm \xi_{k}^v,\xi_{\ell}^p),  \quad 1 \leq k \leq 2N, \quad 1 \leq \ell \leq N,
    \end{aligned}
  \end{equation}
  are independent of the parameters $\param$, and can be precomputed and stored. This is the 
  parameter-independent \emph{offline stage}. In the \emph{online stage},
  for any $\param \in \mathcal{D}$ we assemble and solve the $3N \times 3N$ system for the
  reduced velocity $\bm u_N \in \mathbb{R}^{2N \times 1}$ and reduced pressure 
  $\bm p_N \in \mathbb{R}^{N \times 1}$ s.t.
  \begin{equation}
    \begin{aligned}
      \left( \sum_{m=1}^{M_a} \Theta_a^m(\param) A^m \right) \bm u_N + \left( \sum_{m=1}^{M_b} \Theta_b^m(\param) B^m \right) \boldsymbol{p}_N &= \boldsymbol{F}(\param) \\
      \sum_{m=1}^{M_b} \Theta_b^m(\param) [B^m]^T \bm u_N &= \boldsymbol{G}(\param).
    \end{aligned}
  \end{equation}  
  In the case of pure Dirichlet boundary conditions we have
  \begin{equation}
	\boldsymbol{F}(\param) = \bm f^F - \sum_{m=1}^{M_a} \Theta_a^m(\param) A^m \tilde{\bm u}_0
  \end{equation}
  where the inhomogeneous Dirichlet condition has been treated by lifting the solution 
  of the homogeneous equation into the space with proper boundary conditions. Therefore 
  also the computation of the right-hand side splits affinely in terms of parametric dependence.
  The assembly and solution of the reduced system can be done independent (and thus very
  efficiently) of the dimensions of the finite element spaces for the velocity $\velospace_h$ 
  and pressure $\prespace_h$ respectively. While the offline stage is more expensive compared 
  to solving the finite element problem, if the parametric PDE is evaluated for sufficiently 
  many $\param$ (typically the cutoff point is between 100-500 PDE evaluations), the inexpensive 
  online stage negates the preliminary costs involved in setting up the reduced basis matrices.

  In our case the viscous transformation tensor $\nu_T$ obtained by free-form deformations
  does not satisfy the affine parameterization assumption (\ref{eq:affine}). This difficulty is
  handled by using the Empirical Interpolation Method (EIM) \cite{BarraultEtAl,MadayEtAl,Rozza}, which
  approximates the nonaffinely parameterized tensor with a suitable affinely parameterized one
  by replacing each term of the viscous transformation tensor with an approximate expansion
  \begin{equation} \label{eq:eim}
    [\nu_T]_{i,j} = \sum_{m=1}^{M_{i,j}} \vartheta^{i,j}_m(\param) \zeta_m^{i,j}(\bx) + \varepsilon^{i,j}(\bx,\param),
  \end{equation}
  where the error terms are chosen to be under some tolerance, 
  $|| \varepsilon^{i,j}(\cdot,\param) ||_{L^{\infty}} < \varepsilon^{EIM}_{tol}$ for all $\param \in \mathcal{D}$.
  For this problem we use $\varepsilon^{EIM}_{tol} =$ 1e-5, which is satisfied with $M_{i,j} \leq 22$
  for $i,j=1,2$. The affine decomposition (\ref{eq:affine}) is then recovered. The tolerance 
  should be chosen small enough so that the introduced additional error term does not dominate 
  the error reduced basis approximation \cite{Nguyen,Rozza}. For an a posteriori error estimate 
  of reduced basis approximation for the Stokes problem we refer at the moment to \cite{RovasThesis}.
  
\section{NUMERICAL EXPERIMENTS}
  
\begin{figure}
    \centering
    \subfigure[Displacement $\eta$ and assumed displacement $\hat{\eta}$ at end of coupling Algorithm 1 for (\ref{eq:minproblem3})]{
      \includegraphics[height=3.7cm]{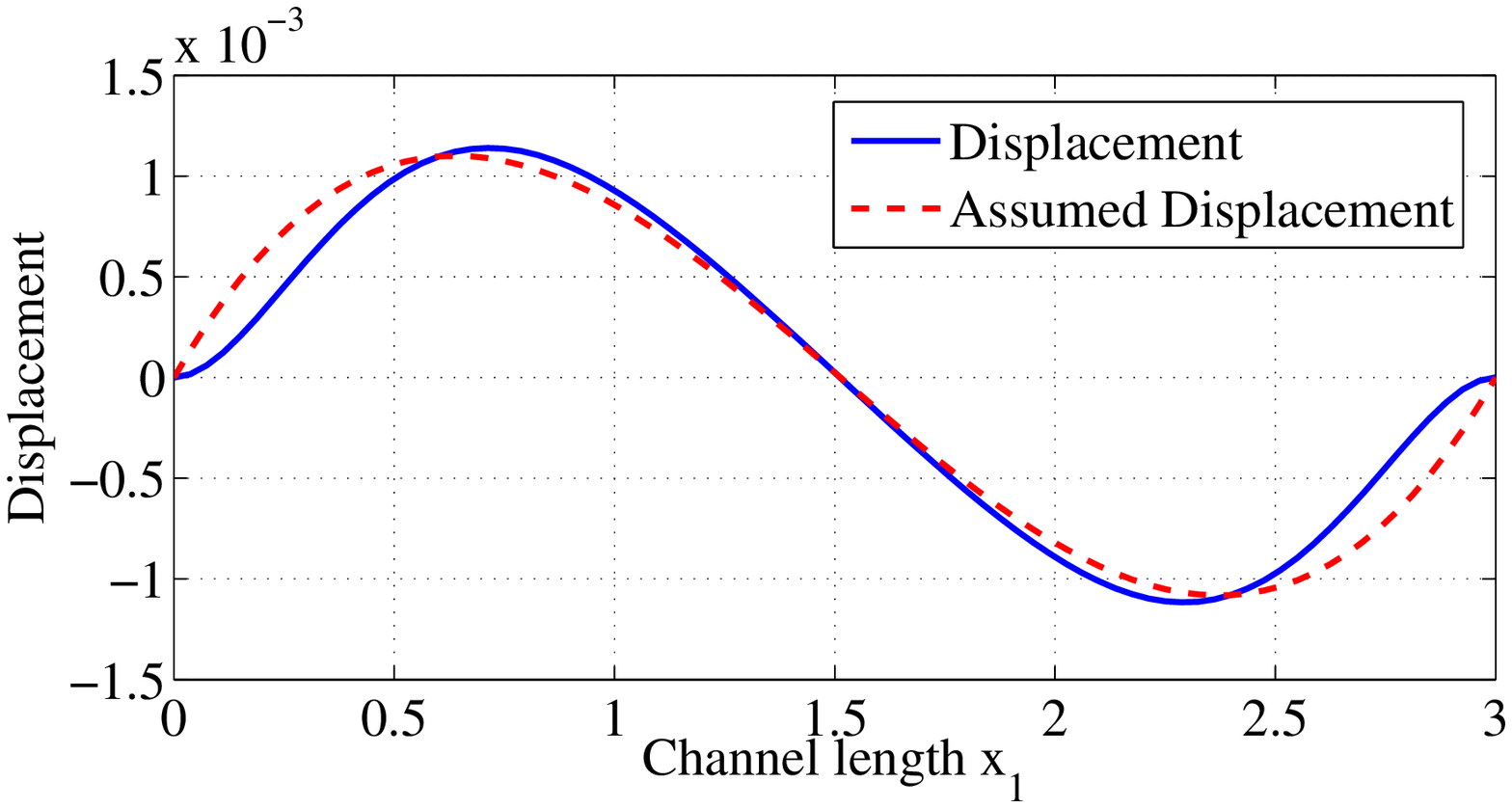}
      \label{fig:subfig1}
    }
    \subfigure[Computational costs of FEM, reduced FEM, and reduced basis approach as number of parametric PDE solutions increases]{
      \includegraphics[height=3.7cm]{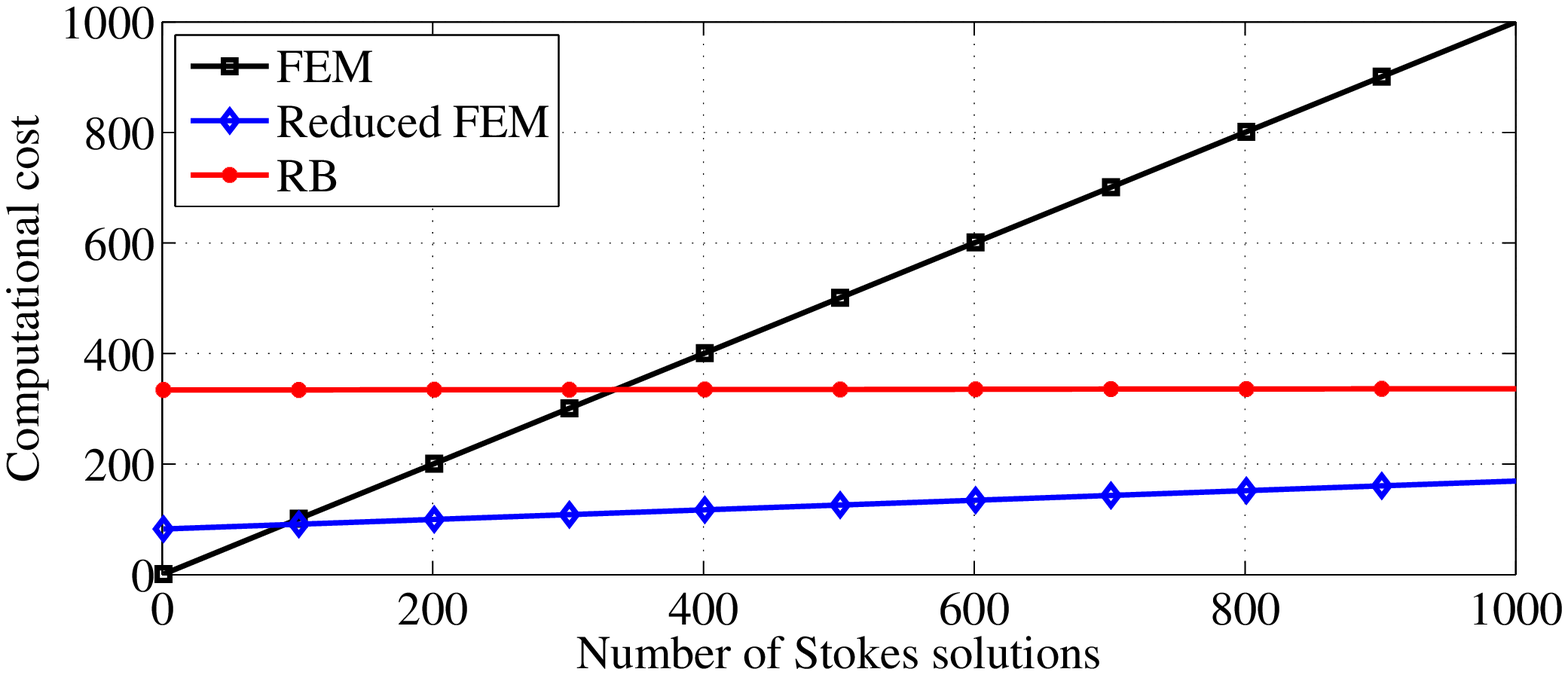}
      \label{fig:subfig2}
    }
    \caption{Result of the fluid-structure algorithm and a comparison of computational costs}
    
  \end{figure}

  Because an analytical solution for the fluid-structure problem described in Sect. \ref{sec:fsi} 
  is not obtainable, we verified our (finite element) computations by comparing them to results 
  presented in \cite{Murea2} for the same test case and parameter values, except with a second-order 
  wall law. For a flexible channel with length 3 cm and half-width 0.5 cm we imposed a Poiseuille 
  flow profile at both the inflow and outflow
  \begin{equation}
    \bm u_0(\bx) = v_0 \begin{bmatrix} 1 - 4 x_2^2 \\ 0 \end{bmatrix}.
  \end{equation}
  The numerical values for the physical parameters were chosen according to \cite{FormaggiaEtAl2} 
  so that the inflow velocity was $v_0 = 30$ cm/s and the blood viscosity $\nu = 0.035$ g/cm$\cdot$s. 
  The volume force on the fluid $\bm f^F$ was taken to be zero. The spring constant was $K = 62.5$ g/s$^2$.
  From the initial guess $\param = \bm 0$ the least-squares optimization algorithm for the strong 
  parametric coupling took 287 iteration steps (Stokes evaluation + solution of least-squares problem
  for the coupling) until the stopping criteria $|| \param^{k+1} - \param^k || < $ 1e-5 was fulfilled. 
  The final value of the cost functional was $J = $ 5e-8. The obtained displacement $\eta$ of the 
  fluid-structure interface versus the assumed displacement $\hat{\eta}$ are shown in Fig. \ref{fig:subfig1}.  

  Once the full FEM solution was verified to be correct, we performed the model reduction of the Stokes problem.   
  We consider two levels of model reduction. In the first level, the problem is parameterized
  with free-form deformations following Sect. \ref{sect:ffd} , the transformation tensors are
  affinely decomposed with the empirical interpolation as in (\ref{eq:eim}), and the parametric
  decomposition (\ref{eq:problem7}) of the system matrices is precomputed.
  Then we can either:
  \begin{enumerate}
    \item Assemble and solve the full FEM Stokes-problem from (\ref{eq:problem7}) using the
    precomputed system matrices for each $\param$. The only reduction is in terms of the 
    parameterized geometry. We call this \emph{reduced FEM}.
    \item Perform the full reduced basis offline stage and construct the reduced system
    matrices in (\ref{eq:problem5}), then solve the reduced basis system for each $\param$.
    Reduction is both in terms of the parameterized geometry and the PDEs. This is the 
    \emph{reduced basis} approach.
  \end{enumerate}
  The computational costs are compared in Fig. \ref{fig:subfig2}. The costs are scaled so that 
  the cost of one full FEM solution equals 1. The full FEM incurs no starting
  overhead, but scales poorly as the number of PDE solutions increases. The reduced FEM has
  moderate starting cost (empirical interpolation + matrix assembly) and scales fairly well.
  The full reduced basis approach is costly in terms of the offline stage, but solving the PDE
  is almost free, resulting in extremely good scalability when thousands of parametric PDE
  evaluations are desired. For the steady FSI problem we don't need thousands of PDE evaluations,
  so that reduced FEM would seem the best option in terms of total
  computational cost for this problem. For the unsteady case we anticipate a need for many more
  PDE solutions during time-stepping so that the added scalability of the full reduced basis
  becomes significant.

\section{CONCLUSIONS}

  A model reduction technique for steady fluid-structure interaction problem of Stokes flow in a
  flexible 2-d channel was proposed. A geometric reduction was performed using free-form deformations
  to reduce the free-boundary problem to a low-dimensional parameter space. Reduced basis
  methods were then used to reduced the complexity of the resulting parametric PDEs.
  
  A least-squares parametric coupling formulation between the fluid and structure was given. The
  approximate balance of stresses across the interface was formulated using the parameterized
  displacement of the interface and solved using nonlinear programming techniques. A low-dimensional 
  parameterization with six parameters coming from a free-form deformation technique was enough
  to obtain approximate coupling.  
  
  Computational costs between the full FEM, a geometrically reduced FEM, and the full reduced basis
  approach were compared. It was observed that the reduced FEM without the overhead of the reduced
  basis reduction was the computationally most attractive choice, but that for unsteady or nonlinear
  problems there could be a need for the better scalability of reduced basis methods as the number
  of PDE evaluations increases.

\section*{ACKNOWLEDGEMENTS}

  We thank Alfio Quarteroni for his comments and suggestions regarding the coupling algorithm.
  The reduced basis computations were performed with the rbMIT toolkit \cite{rbmit}.
  Andrea Manzoni contributed code for the numerical simulations. A previous version of this article
  appeared in: R.A.E.~M{\"{a}}kinen, P.~Neittaanm{\"{a}}ki, T.~Tuovinen, K.~Valpe
  (Eds.) Proceedings of the 10th Finnish Mechanics Days, December 3-4, Jyv{\"{a}}skyl{\"{a}}, Finland, 2010.
  Published with permission.

  \bibliography{my}
  \bibliographystyle{plain}

\end{document}

%% file: fsi.pstex_t
\begin{picture}(0,0)%
\includegraphics{fsi.pstex}%
\end{picture}%
\setlength{\unitlength}{3947sp}%
\begingroup\makeatletter\ifx\SetFigFont\undefined%
\gdef\SetFigFont#1#2#3#4#5{%
  \reset@font\fontsize{#1}{#2pt}%
  \fontfamily{#3}\fontseries{#4}\fontshape{#5}%
  \selectfont}%
\fi\endgroup%
\begin{picture}(6180,1908)(2311,-3673)
\put(4501,-2386){\makebox(0,0)[lb]{\smash{{\SetFigFont{12}{14.4}{\rmdefault}{\mddefault}{\updefault}$\eta(x_1)$}}}}
\put(2326,-3136){\makebox(0,0)[lb]{\smash{{\SetFigFont{12}{14.4}{\rmdefault}{\mddefault}{\updefault}$\Gamma_{in}$}}}}
\put(8476,-3136){\makebox(0,0)[lb]{\smash{{\SetFigFont{12}{14.4}{\rmdefault}{\mddefault}{\updefault}$\Gamma_{out}$}}}}
\put(4051,-3286){\makebox(0,0)[lb]{\smash{{\SetFigFont{12}{14.4}{\rmdefault}{\mddefault}{\updefault}$\Omega(\eta)$}}}}
\put(5926,-2986){\makebox(0,0)[lb]{\smash{{\SetFigFont{12}{14.4}{\rmdefault}{\mddefault}{\updefault}$-\triangle \boldsymbol{u} + \nabla p = 0$}}}}
\put(6489,-3278){\makebox(0,0)[lb]{\smash{{\SetFigFont{12}{14.4}{\rmdefault}{\mddefault}{\updefault}$\nabla \cdot \boldsymbol{u} = 0$}}}}
\put(4951,-1936){\makebox(0,0)[lb]{\smash{{\SetFigFont{12}{14.4}{\rmdefault}{\mddefault}{\updefault}$\Sigma(\eta)$}}}}
\end{picture}%